\documentclass[12pt,reqno]{amsart}
\newcommand\ignore[1]{}
\input amssym.def
\usepackage[raiselinks=false,colorlinks=true,citecolor=blue,urlcolor=blue,linkcolor=blue,bookmarksopen=true,pdftex]{hyperref}
\usepackage{booktabs,multirow,lscape,longtable} % for tables 
\usepackage{mathtools} 
\usepackage{tikz} 
\numberwithin{equation}{section}
\numberwithin{equation}{subsection} 
\newtheorem{theorem}{Theorem}[section]

\newtheorem{remark}[theorem]{Remark}

\begin{document} 
\baselineskip=15.5pt

\title
[Irreducible unitary representations with non-zero $(\frak{g}, K)$-cohomology]{Irreducible unitary representations with non-zero relative Lie algebra cohomology of a Lie group of 
type $\frak{f}_{4(4)}$}  
\author{Pampa Paul }
\address{Department of Mathematics, Presidency University, 86/1 College Street, Kolkata 700073, India}
\email{pampa.maths@presiuniv.ac.in }
\subjclass[2020]{22E46, 17B10, 17B20, 17B22, 17B25, 17B56.  \\ 
Keywords and phrases: Lie group, Lie algebra, Dynkin diagram, $\theta$-stable parabolic subalgebra, cohomological induction.}

\thispagestyle{empty}
\date{}

\begin{abstract}
In this article, we have determined the irreducible unitary representations with non-zero relative Lie algebra cohomology and Poincar\'{e} polynomials of 
cohomologies of these representations for a connected Lie group $G$ with Lie algebra $\frak{f}_{4(4)}.$ We have also determined a necessary and sufficient condition 
for these representations to be discrete series representations and identified the discrete series representations and Borel-de Siebenthal discrete series representations 
among the irreducible unitary representations of $G$ with non-zero relative Lie algebra cohomology. 
\end{abstract}
\maketitle

\noindent 
\section{Introduction} 

Let $G$ be a connected semisimple Lie group with finite centre, and $K$ be a maximal compact subgroup of $G$ with Cartan involution $\theta.$ The differential of $\theta$ at identity is denoted 
by the same notation $\theta.$ Let $\frak{g}_0=$Lie$(G),\frak{k}_0=$Lie$(K),$ and $\frak{g}_0=\frak{k}_0 \oplus \frak{p}_0$ be the Cartan decomposition corresponding to $\theta.$ Let 
$\frak{g}=\frak{g}_0^\mathbb{C},\frak{k}=\frak{k}_0^\mathbb{C}, \frak{p}=\frak{p}_0^\mathbb{C}.$ Assume that $\frak{h}_0$ be a $\theta$-stable 
fundamental Cartan subalgebra of $\frak{g}_0,$ and $\frak{h}=\frak{h}_0^\mathbb{C}.$ Corresponding to a $\theta$-stable parabolic subalgebra $\frak{q}$ of $\frak{g}_0$ containing $\frak{h}_0,$ 
and a linear function $\lambda$ on $\frak{h}$ in a certain good range, there is a 
cohomologically induced module whose completion $A_\frak{q}(\lambda)$ is an irreducible unitary representation of $G$ with infinitesimal character $\chi_\lambda.$ 
 If rank$(G)=$ rank$(K),$ any discrete series representation of $G$ is an $A_\frak{q}(\lambda).$ Among these, the $A_\frak{q}(0)$s are of particular interest because these are  
irreducible unitary representations with non-zero $(\frak{g}, K)$-cohomology, and any irreducible unitary representation with non-zero $(\frak{g}, K)$-cohomology has this form. 
We denote $A_\frak{q}(0)$ by $A_\frak{q}.$ 

A $\theta$-stable parabolic subalgebra of $\frak{g}_0$ is a parabolic subalgebra $\frak{q}$ of $\frak{g}$ such that 
(a) $\theta(\frak{q}) = \frak{q}$, and (b) $\bar{\frak{q}} \cap \frak{q}= \frak{l}$ is a Levi subalgebra of $\frak{q}$; 
where $\bar{\ }$ denotes the conjugation of $\frak{g}$ with respect to $\frak{g}_0$. If $\frak{q}$ is a $\theta$-stable parabolic subalgebra, then by (b), $\frak{l}$ is the 
complexification of a real subalgebra $\frak{l}_0$ of $\frak{g}_0$ and $\frak{l}_0$ is $\theta$-stable. Let $\frak{u}$ be the nilradical of $\frak{q}$ so that $\frak{q} = \frak{l} \oplus \frak{u}.$ 
Then $\frak{u}$ is $\theta$-stable and so $\frak{u} = (\frak{u} \cap \frak{k}) \oplus (\frak{u} \cap \frak{p}).$ Now associated with a $\theta$-stable parabolic subalgebra $\frak{q}$, 
we have an irreducible unitary representation $A_\frak{q}$ whose $(\frak{g},K)$-module $A_{\frak{q},K}$ is given by $\mathcal{R}^S _\frak{q} (\mathbb{C}),$ where $S = \textrm{dim}
(\frak{u} \cap \frak{k}).$ The irreducible unitary representation $A_\frak{q}$ has trivial infinitesimal character. The Levi subgroup $L = \{g \in G : \textrm{Ad}(g) (\frak{q}) = 
\frak{q} \}$ is a connected reductive Lie subgroup of $G$ with Lie algebra $\frak{l}_0.$ As $\theta(\frak{l}_0) = \frak{l}_0, L \cap K$ is a maximal 
compact subgroup of $L$. We have  
\[ H^r (\frak{g}, K; A_{\frak{q}, K}) \cong H^{r-R(\frak{q})} (\frak{l}, L\cap K ; \mathbb{C}), \]
where $R(\frak{q}) := \textrm{dim}(\frak{u} \cap \frak{p})$. Let $Y_\frak{q}$ denote the compact dual of the 
Riemannian globally symmetric space $L/{L\cap K}$. Then $H^r (\frak{l}, L\cap K ; \mathbb{C}) \cong 
H^r (Y_\frak{q} ; \mathbb{C})$. And hence 
\[  H^r (\frak{g}, K; A_{\frak{q}, K}) \cong H^{r-R(\frak{q})} (Y_\frak{q} ; \mathbb{C}).\]
If $P_\frak{q}(t)$ denotes the Poincar\'{e} polynomial of $ H^* (\frak{g}, K; A_{\frak{q}, K})$, then 
\[ P_\frak{q}(t) = t^{R(\frak{q})} P(Y_\frak{q} , t). \]
 Conversely, if $\pi$ is an irreducible unitary represention of $G$ with non-zero $(\frak{g},K)$-cohomology,  then $\pi $ is unitarily equivalent 
to $A_\frak{q}$ for some $\theta$-stable parabolic subalgebra $\frak{q}$ of $\frak{g}_0$ \cite[Th. 4.1]{voganz}. 

If rank$(G) =$ rank$(K)$ and $\frak{q}$ is a $\theta$-stable Borel subalgebra, then $A_\frak{q}$ is a 
discrete series representation of $G$ with trivial infinitesimal character. In this case, $R(\frak{q})= \frac{1}{2} \textrm{ dim}(G/K)$,  
$\frak{l}_0 = \frak{t}_0$ and hence 
\[H^r (\frak{g}, K; A_{\frak{q}, K}) =
\begin{cases}
0 & \textrm{if } r \neq R(\frak{q}), \\
\mathbb{C} & \textrm{if } r = R(\frak{q}). 
\end{cases}
\]
If we take $\frak{q} = \frak{g}$, then $L=G$ and $A_\frak{q} = \mathbb{C}$, the trivial representation of $G$. 

If $\frak{q}$ is a $\theta$-stable parabolic subalgebra, then so is $Ad(k)(\frak{q})$ and 
$A_\frak{q}$ is unitarily equivalent to $A_{Ad(k)(\frak{q})}$ for all $k \in K.$ Now the $\theta$-stable real form $\frak{l}_0$ of the Levi subalgebra $\frak{l}$ contains a 
maximal abelian subalgebra $\frak{t}_0$ of $\frak{k}_0$. Then $\frak{z}_{\frak{g}_0} (\frak{t}_0)$ is a $\theta$-stable fundamental Cartan subalgebra of $\frak{g}_0.$
Since any two fundamental Cartan subalgebras of $\frak{g}_0$ are $Ad(k)$-conjugate for some $k \in K,$ we may assume that $\frak{z}_{\frak{g}_0} (\frak{t}_0) = \frak{h}_0.$
Let $\frak{t}=\frak{t}_0^\mathbb{C}.$ Note that $\frak{t},\frak{h}$ are Cartan subalgebras of $\frak{k},\frak{g}$ respectively and $\frak{h} \subset \frak{q}.$ Let 
$\Delta = \Delta(\frak{g}, \frak{h})$ be the set of all non-zero roots of $\frak{g}$ relative to the Cartan subalgebra $\frak{h},$ and $\Delta_\frak{k} = \Delta(\frak{k},\frak{t}).$ 
Choose a system of positive roots $\Delta_\frak{k}^+$ in $\Delta_\frak{k}.$ If $x \in i\frak{t}_0$ be such that $\alpha(x) \ge 0$ for all $\alpha \in \Delta_\frak{k}^+$, then 
$\frak{q}_x= \frak{h}\oplus \sum_{\alpha \in \Delta,\alpha(x) \ge 0} \frak{g}^\alpha$ 
is a $\theta$-stable parabolic subalgebra of $\frak{g}_0$, $\frak{l}_x= \frak{h}\oplus \sum_{\alpha \in \Delta,\alpha(x)= 0} \frak{g}^\alpha$ is the 
Levi subalgebra of $\frak{q}_x,$ and $\frak{u}_x= \sum_{\alpha \in \Delta,\alpha(x) > 0} \frak{g}^\alpha$ is the nilradical of $\frak{q}_x,$ where 
$\frak{g}^\alpha$ is the root subspace of $\frak{g}$ corresponding to the non-zero root $\alpha \in \Delta.$ If $\frak{q}$ is a 
$\theta$-stable parabolic subalgebra of $\frak{g}_0$, there exists $k \in K$ such that $Ad(k)(\frak{q})=\frak{q}_x.$ Let $\mathcal{Q}$ denote the set of all $\theta$-stable 
parabolic subalgebras of the form $\frak{q}_x,$ where $x \in i\frak{t}_0$ is such that $\alpha(x) \ge 0$ for all $\alpha \in \Delta_\frak{k}^+.$ If $\frak{q}$ is a 
$\theta$-stable parabolic subalgebra, then $\frak{u} \cap \frak{p}$ is a $\frak{t}$-module. Let $\Delta(\frak{u} \cap \frak{p})$ denote the set of all non-zero weights of 
$\frak{u} \cap \frak{p},$ and $\delta (\frak{u} \cap \frak{p})$ denote $1/2$ of the sum of elements in $\Delta (\frak{u} \cap \frak{p}).$ 
For $\frak{q}, \frak{q}' \in \mathcal{Q}, A_\frak{q}$ is unitarily equivalent to $A_{\frak{q}'}$ {\it if and only if} $\Delta(\frak{u} \cap \frak{p})=
\Delta(\frak{u'} \cap \frak{p})$ \cite{riba}. 

The associated $(\frak{g},K)$-module $A_{\frak{q}, K}$ of $A_\frak{q}$ contains an 
irreducible $K$-submodule of highest weight (with respect to $\Delta ^+ _\frak{k}$) 
$2 \delta (\frak{u} \cap \frak{p}) = 
\sum_{\beta \in \Delta (\frak{u} \cap \frak{p})} \beta $ and 
it occurs with multiplicity one in $A_{\frak{q}, K}$. Any other irreducible $K$-module that occurs in $A_{\frak{q}, K}$ has highest weight 
of the form $2 \delta (\frak{u} \cap \frak{p}) + 
\sum_{\gamma \in \Delta (\frak{u} \cap \frak{p})} n_\gamma \gamma,$
with $n_\gamma$ a non-negative integer \cite[Th. 2.5]{voganz}.  
The $(\frak{g}, K)$-modules $A_{\frak{q} , K}$ were first constructed, in general, by Parthasarathy \cite{parthasarathy1}. 
Vogan and Zuckerman \cite{voganz} gave a construction of the $(\frak{g}, K)$-modules $A_{\frak{q} , K}$ via cohomological induction and 
Vogan \cite{vogan} proved that these are unitarizable. In this article, we have determined a necessary and sufficient condition for an $A_\frak{q}$ 
to be a discrete series representation with trivial infinitesimal character, in terms of $\Delta(\frak{u} \cap \frak{p}).$ The result is given below and the 
proof is given in \S \ref{discrete}.  

\begin{theorem}\label{th1}
Assume that rank$(G)=$ rank$(K).$ Then a representation $A_\frak{q}$ is a discrete series representation with trivial infinitesimal character 
{\it if and only if} $\Delta_\frak{k}^+ \cup \Delta(\frak{u}\cap \frak{p})$ is a positive root system of $\Delta.$ If $G/K$ is a Hermitian symmetric space, then 
$A_\frak{q}$ is a holomorphic discrete series representation with trivial infinitesimal character {\it if and only if} 
$\Delta_\frak{k}^+ \cup \Delta(\frak{u}\cap \frak{p})$ is a Borel-de Siebenthal positive root system of $\Delta.$ 
If $G/K$ is not Hermitian symmetric, then $A_\frak{q}$ is a Borel-de Siebenthal discrete series representation with trivial infinitesimal character {\it if and only if} 
$\Delta_\frak{k}^+ \cup \Delta(\frak{u}\cap \frak{p})$ is a Borel-de Siebenthal positive root system of $\Delta.$ 
\end{theorem}

Millson and Raghunathan \cite{mira} have determined a technique to identify certain $A_\frak{q}$s as automorphic representations 
by constructing geometric cycles and using Matsushima's isomorphism \cite{matsushima}. Later this technique was used by various authors. 
Collingwood \cite{collingwood} has determined the cohomologies of the representations $A_\frak{q}$ of $Sp(n,1),$ and the real rank one real form of $F_4.$ 
Li and Schwermer \cite{lisch} have computed the cohomologies of these representations for the connected non-compact real Lie group of 
type $G_2.$ Mondal and Sankaran \cite{mondal-sankaran2} have determined the representations $A_\frak{q}$ of Hodge type $(p,p)$ for small values of $p,$ 
when $G/K$ is an irreducible Hermitian symmetric space. 
If $G$ is a complex simple Lie group, or $G=SO_0(2,m)(m \in \mathbb{N}),$ 
the number of unitary equivalence classes of the representations $A_\frak{q},$ and Poincar\'{e} polynomials of cohomologies of these representations have been 
determined in \cite{paul}, \cite{paul1}. In this article, we have determined the number of unitary equivalence classes of the representations $A_\frak{q},$ 
the number of unitary equivalence classes of discrete series representations with trivial infinitesimal character, the number of unitary equivalence classes of 
Borel-de Siebenthal discrete series representations with trivial infinitesimal character, and Poincar\'{e} polynomials of cohomologies of these representations,  
when $G$ is a connected Lie group with Lie algebra $f_{4(4)}.$ The process of determining the inequivalent $A_\frak{q}$s in this article is different from \cite{paul} and 
\cite{paul1}, due to the complexity of relation between the non-compact roots. 
Here the Lie algebra $\frak{f}_{4(4)}$ is the non-compact real form of the exceptional complex simple Lie algebra $\frak{f}_4$ such that the 
maximal compactly imbedded subalgebra of $\frak{f}_{4(4)}$ is $\frak{sp}(3) \oplus \frak{su}(2).$ Also rank$(\frak{f}_{4(4)})=$ rank$(\frak{sp}(3) \oplus \frak{su}(2)),$ 
$\frak{sp}(3) \oplus \frak{su}(2)$ is semisimple, and $\frak{f}_{4(4)}$ has real rank $4.$ 

\begin{theorem}\label{th2}
Let $G$ be a connected Lie group with Lie algebra $\frak{f}_{4(4)}.$ Then the number of unitary equivalence classes of irreducible unitary representations of $G$ with 
non-zero $(\frak{g}, K)$-cohomology is $46.$ Among these, the unitary equivalence classes of discrete series representations with trivial infinitesimal character are 
$12$ in number and exactly one of them is the unitary equivalence class of Borel-de Siebenthal discrete series representations with trivial infinitesimal character. The 
Poincar\'{e} polynomials of cohomologies of these representations are listed in the Table \ref{f-table}. 
\end{theorem}

The proof is given in \S \ref{f}. 

\noindent 
\section{Discrete series representations with trivial infinitesimal character}\label{discrete}

We follow the notations from the previous section. In addition, assume that rank$(G)=$ rank$(K).$ Let $W_\frak{g}$ (respectively $W_\frak{k}$) denote 
the Weyl group of $\frak{g}$ (respectively, $\frak{k}$) relative to the Cartan subalgebra $\frak{h}.$ A non-singular linear function $\lambda$ on $i\frak{t}_0$ relative to 
$\Delta$ defines uniquely a positive root system $\Delta_\lambda^+$ of $\Delta.$ Define $\delta_\frak{g} = \frac{1}{2}\sum_{\alpha \in \Delta_\lambda^+}\alpha, \delta_\frak{k} = \frac{1}{2} 
\sum_{\alpha \in \Delta_\lambda^+ \cap \Delta_\frak{k}}\alpha.$ If $\lambda +\delta_\frak{g}$ is 
analytically integral(that is, $\lambda +\delta_\frak{g}$ is the differential of a Lie group homomorphism on the Cartan subgroup of $G$ corresponding to $\frak{t}_0$), then there exists a 
discrete series representation $\pi_\lambda$ with infinitesimal character $\chi_\lambda;$ the associated $(\frak{g}, K)$-module $\pi_{\lambda,K}$ contains an 
irreducible $K$-submodule with highest weight $\Lambda=\lambda+\delta_\frak{g}-2\delta_\frak{k};$
it occurs with multiplicity one in $\pi_{\lambda, K}.$ Any other irreducible $K$-module that occurs in $\pi_{\lambda, K}$ has highest weight 
of the form $\lambda + \sum_{\alpha \in \Delta_\lambda^+} n_\alpha \alpha,$ with $n_\alpha$ a non-negative integer; and 
two such representations $\pi_\lambda,$ and $\pi_{\lambda'}$ are unitarily equivalent {\it iff} $\lambda = w\lambda'$ for some $w \in W_\frak{k}$ \cite[Th. 9.20, Ch. IX]{knapp}. 
This $\lambda$ is called the {\it Harish-Chandra parameter}, $\Lambda$ is called the {\it Blattner parameter} of the discrete series representation $\pi_\lambda,$ 
and the positive root system $\Delta_\lambda^+$ is called the {\it Harish-Chandra root order} corresponding to $\lambda.$
Upto unitary equivalence, these are all discrete series representations of $G$ \cite[Th. 12.21, Ch. XII]{knapp}. To get non-equivalent discrete series 
representations, we assume that the Harish-Chandra root order $\Delta_\lambda^+$ corresponding to $\lambda$ contains $\Delta_\frak{k}^+$ 
so that $\Delta_\lambda^+ \cap \Delta_\frak{k} = \Delta_\frak{k}^+.$ 

Since the Cartan subalgebra $\frak{h} \subset \frak{k}$ and each root subspace $\frak{g}^\alpha (\alpha \in \Delta(\frak{g}, \frak{h}))$ is one-dimensional, either 
$\frak{g}^\alpha \subset \frak{k}$ of $\frak{g}^\alpha \subset \frak{p}.$ A root $\alpha$ is compact if $\frak{g}^\alpha \subset \frak{k},$ and is non-compact if 
$\frak{g}^\alpha \subset \frak{p}.$ If $G$ is simple, a {\it Borel-de Siebenthal positive root system} is a positive root system of $\Delta(\frak{g}, \frak{h})$ such that it contains 
exactly one non-compact simple root $\nu$ and the coefficient of $\nu$ in the highest root when expressed as the sum of simple roots, is $1$ if $G/K$ is 
Hermitian symmetric, and is $2$ if $G/K$ is not Hermitian symmetric. If $G$ is semisimple, then a Borel-de Siebenthal positive root system is a positive root system of 
$\Delta(\frak{g}, \frak{h})$ such that it is a union of Borel-de Siebenthal positive root systems of each simple ideal of $\frak{g}.$ 
If $G/K$ is Hermitian symmetric, then $\pi_\lambda$ is a 
holomorphic discrete series representation {\it iff} the Harish-Chandra root order corresponding to $\lambda$ is a Borel-de Siebenthal positive root system \cite{hc1}. 
If $G/K$ is not Hermitian symmetric, then $\pi_\lambda$ is a Borel-de Siebenthal discrete series representation 
(defined in \cite{ow} analogous to holomorphic discrete series representations) {\it iff} the Harish-Chandra root order $\Delta_\lambda^+$ is a 
Borel-de Siebenthal positive root system \cite{ow}. Note that the infinitesimal character $\chi_\lambda$ is the character of the Verma module of 
$\frak{g}$ with highest weight $\lambda-\delta_\frak{g}.$ Thus $\chi_\lambda$ is trivial {\it iff} $\lambda = \delta_\frak{g}$ 
{\it iff} $\Lambda = \sum_{\beta \in (\Delta_\lambda^+ \setminus \Delta_\frak{k}^+)} \beta.$ Now the proof of Th. \ref{th1} is given below. 

{\bf Proof of Theorem \ref{th1}: } 
The representation $A_\frak{q}$ is a discrete series representation with trivial infinitesimal character {\it if and only if} $A_\frak{q}$ is unitarily equivalent to $A_\frak{b}$ for some 
Borel subalgebra $\frak{b}$ of $\frak{g}$ containing the Borel subalgebra  $\frak{b}_\frak{k}=\frak{h} \oplus \sum_{\alpha \in \Delta_\frak{k}^+} \frak{g}^\alpha$ of $\frak{k}.$ And in this case, the 
Blattner parameter $\Lambda$ of the discrete series representation  $A_\frak{q}$ is $\sum_{\beta \in \Delta(\frak{u} \cap \frak{p})} \beta,$ and so the Harish-Chandra root 
order is $\Delta_\frak{k}^+ \cup \Delta(\frak{u} \cap \frak{p}).$ Thus $\Delta_\frak{k}^+ \cup \Delta(\frak{u} \cap \frak{p})$ is a positive root system of $\Delta.$ To prove the converse 
part, let $\frak{b} = \frak{h} \oplus \sum_{\beta \in \Delta_\frak{k}^+\cup \Delta(\frak{u}\cap \frak{p})} \frak{g}^\beta.$ Then $\frak{b}$ is a Borel subalgebra of $\frak{g}$ containing 
$\frak{b}_\frak{k},$ and $\Delta(\frak{u}_\frak{b} \cap \frak{p})=\Delta(\frak{u} \cap \frak{p}), \frak{u}_\frak{b}$ is the nilradical of $\frak{b}.$ Thus $A_\frak{q}$ is unitarily 
equivalent to $A_\frak{b},$ which is a discrete series representation. This proves the first part. The second and the third part follow from the first part.

\noindent
\section{Irreducible unitary representations with non-zero $(\frak{g}, K)$-cohomology of a Lie group of type $\frak{f}_{4(4)}$}\label{f}  

Let $G$ be a connected Lie group with Lie algebra $\frak{g}_0= \frak{f}_{4(4)},$ and $\frak{g}_0 = \frak{k}_0 \oplus \frak{p}_0$ be a Cartan decomposition of $\frak{g}_0$ with 
the corresponding Cartan involution $\theta.$ Then $\frak{k}_0 \cong \frak{sp}(3) \oplus \frak{su}(2).$ Let $K$ be the connected Lie subgroup of $G$ with Lie algebra $\frak{k}_0.$ Since the 
centre of $G$ is finite, $K$ is a maximal compact subgroup of $G.$ The Lie group automorphism of $G$ whose differential at the identity is $\theta,$ is denoted by the same notation $\theta.$ 
Note that $K$ is semisimple, rank$(G)=$rank$(K),$ and $G/K$ is an irreducible Riemannian globally symmetric space of non-compact type. 

Let $\frak{g}=\frak{g}_0^\mathbb{C}; \frak{k}=\frak{k}_0^\mathbb{C}, \frak{p}=\frak{p}_0^\mathbb{C} \subset \frak{g}.$ Then $\frak{g}$ is the simple Lie algebra $\frak{f}_4.$ Let 
$\frak{t}_0$ be a maximal abelian subspace of $\frak{k}_0$. Then $\frak{h}= \frak{t}_0^\mathbb{C} \subset \frak{k}$ is a Cartan subalgebra of $\frak{k}$ as well as of $\frak{g}.$ Let 
$\Delta = \Delta (\frak{g}, \frak{h})$ be the set of non-zero roots of $\frak{g}$ relative to the Cartan subalgebra $\frak{h},$ and similarly $\Delta_\frak{k}=\Delta(\frak{k}, \frak{h}).$ 
Note that $\Delta_\frak{k} \subset \Delta,$ and $\Delta_\frak{k}$ is the set of all compact roots in $\Delta.$ Let $\Delta_n = \Delta \setminus \Delta_\frak{k},$ and it is the set of all 
non-compact roots in $\Delta.$ Let $\Delta^+$ is a Borel-de Siebenthal positive root system of $\Delta$ with a unique non-compact simple root, say $\nu.$ Then $n_\nu (\alpha) = 0 
\textrm{ or } \pm 2,$ if $\alpha \in \Delta_\frak{k},$ and $n_\nu(\alpha)=\pm 1,$ if $\alpha \in \Delta_n;$ where $n_\nu (\alpha)$ is the coefficient of $\nu$ in the decomposition of $\alpha$ 
in terms of the simple roots of $\Delta^+.$ If $\phi_1, \phi_2, \phi_3, \phi_4$ are the simple roots of $\Delta^+$ with $\nu = \phi_1,$ and $\delta$ is the highest root, then 
$\delta = 2\phi_1 +3\phi_2 +4\phi_3 +2\phi_4.$ The extended Dynkin diagram of $\frak{g}=\frak{f}_4$ is given by

\begin{center} 
\begin{tikzpicture}

\draw (0,0) circle [radius = 0.1]; 
\filldraw [black] (1,0) circle [radius = 0.1]; 
\draw (2,0) circle [radius = 0.1]; 
\draw (3,0) circle [radius = 0.1]; 
\draw (4,0) circle [radius = 0.1]; 
\node [above] at (0,0.05) {$-\delta$}; 
\node [above] at (1.05,0.05) {$\phi_1$}; 
\node [above] at (2.05,0.05) {$\phi_2$}; 
 \node [above] at (3.05,0.05) {$\phi_3$}; 
\node [above] at (4.05,0.05) {$\phi_4$}; 
\node [left] at (-0.3,0) {$\frak{f}_4^{(1)} :$}; 
\draw (0.1,0) -- (0.9,0); 
\draw (1.1,0) -- (1.9,0);
\draw (2.9,0) -- (2.8,0.1); 
\draw (2.9,0) -- (2.8,-0.1); 
\draw (2.1,0.025) -- (2.85,0.025); 
\draw (2.1,-0.025) -- (2.85,-0.025); 
\draw (3.1,0) -- (3.9,0); 

 \end{tikzpicture} 
\end{center}

In the diagrams of this article the non-compact roots are represented by black vertices. Note that the Dynkin diagram of $\frak{k}$ is the subdiagram of the extended Dynkin diagram of 
$\frak{f}_4$ with the vertex $\phi_1$ deleted. 

Let $\Delta_\frak{k}^+ =\Delta^+ \cap \Delta_\frak{k}, \Delta_n^+ = \Delta^+ \cap \Delta_n.$ Then $\Delta_\frak{k}^+ = \Delta_0^+ \cup \Delta_2, \Delta_n^+ = \Delta_1;$ 
where $\Delta_i = \{ \alpha \in \Delta : n_\nu (\alpha) = i \}$ for $i=0,\pm1,\pm2,$ and $\Delta_0^+ = \Delta^+ \cap \Delta_0.$ 
Note that $\Delta_n^+= \{\phi_1, \phi_1 +\phi_2, \phi_1+\phi_2+\phi_3, \phi_1+\phi_2+2\phi_3, 
\phi_1+\phi_2+\phi_3+\phi_4, \phi_1+2\phi_2+2\phi_3, \phi_1+\phi_2+2\phi_3+\phi_4, \phi_1+2\phi_2+2\phi_3+\phi_4, \phi_1+\phi_2+2\phi_3+2\phi_4, \phi_1+2\phi_2+3\phi_3+\phi_4, 
\phi_1+2\phi_2+2\phi_3+2\phi_4, \phi_1+2\phi_2+3\phi_3+2\phi_4, \phi_1+2\phi_2+4\phi_3+2\phi_4, \phi_1+3\phi_2+4\phi_3+2\phi_4 \},$ and $\Delta_2 = \{ 
\delta=2\phi_1+3\phi_2+4\phi_3+2\phi_4\}.$ So the set of all simple roots in $\Delta_k^+$ is given by $\{\phi_2, \phi_3 , \phi_4, \delta\}.$

\begin{figure}[!h]
\begin{center} 
\begin{tikzpicture}

\filldraw [black] (0,0) circle [radius = 0.1]; 
\filldraw [black] (0,-1.5) circle [radius = 0.1]; 
\filldraw [black] (0,-3) circle [radius = 0.1]; 
\filldraw [black] (-1,-4) circle [radius = 0.1]; 
\filldraw [black] (1,-4) circle [radius = 0.1]; 
\filldraw [black] (-2,-5) circle [radius = 0.1]; 
\filldraw [black] (0,-5) circle [radius = 0.1]; 
\filldraw [black] (-1,-6) circle [radius = 0.1]; 
\filldraw [black] (1,-6) circle [radius = 0.1]; 
\filldraw [black] (-2,-7) circle [radius = 0.1]; 
\filldraw [black] (0,-7) circle [radius = 0.1]; 
\filldraw [black] (-1,-8) circle [radius = 0.1]; 
\filldraw [black] (-1,-9.5) circle [radius = 0.1]; 
\filldraw [black] (-1,-11) circle [radius = 0.1]; 

\node [right] at (0,0) {$\phi_1$}; 
\node [right] at (0,-1.5) {$\phi_1+\phi_2$}; 
\node [right] at (0,-3) {$\phi_1+\phi_2+\phi_3$}; 
\node [left] at (-1,-4) {$\phi_1+\phi_2+2\phi_3$}; 
\node [right] at (1,-4) {$\phi_1+\phi_2+\phi_3 +\phi_4$}; 
\node [left] at (-2,-5) {$\phi_1+2\phi_2+2\phi_3$}; 
\node [right] at (0,-5) {$\phi_1+\phi_2+2\phi_3 +\phi_4$}; 
\node [left] at (-1,-6) {$\phi_1+2\phi_2+2\phi_3+\phi_4$}; 
\node [right] at (1,-6) {$\phi_1+\phi_2+2\phi_3 +2\phi_4$}; 
\node [left] at (-2,-7) {$\phi_1+2\phi_2+3\phi_3+\phi_4$}; 
\node [right] at (0,-7) {$\phi_1+2\phi_2+2\phi_3 +2\phi_4$}; 
\node [right] at (-1,-8) {$\phi_1+2\phi_2+3\phi_3 +2\phi_4$}; 
\node [right] at (-1,-9.5) {$\phi_1+2\phi_2+4\phi_3 +2\phi_4$};
\node [right] at (-1,-11) {$\phi_1+3\phi_2+4\phi_3 +2\phi_4$};

\node [left] at (-6,-5.5) {$\Delta_n^+ :$}; 

\node [right] at (-0.1,-0.7) {$\phi_2$}; 
\node [right] at (-0.1,-2.2) {$\phi_3$};
\node [left] at (-0.4,-3.4) {$\phi_3$};
\node [right] at (0.4,-3.4) {$\phi_4$};
\node [right] at (-0.7,-4.3) {$\phi_4$};
\node [left] at (0.7,-4.3) {$\phi_3$};
\node [left] at (-1.4,-4.4) {$\phi_2$};
\node [right] at (-1.7,-5.3) {$\phi_4$};
\node [left] at (-0.3,-5.3) {$\phi_2$};
\node [right] at (0.4,-5.4) {$\phi_4$};
\node [left] at (-1.4,-6.4) {$\phi_3$};
\node [right] at (-0.7,-6.3) {$\phi_4$};
\node [left] at (0.7,-6.3) {$\phi_2$};
\node [right] at (-1.7,-7.3) {$\phi_4$};
\node [left] at (-0.3,-7.3) {$\phi_3$};
\node [right] at (-1.1,-8.7) {$\phi_3$};
\node [right] at (-1.1,-10.2) {$\phi_2$};

\draw (0,-0.1) -- (0,-1.4); 
\draw (0,-1.6) -- (0,-2.9);
\draw (0,-3.1) -- (-1,-3.9);
\draw (0,-3.1) -- (1,-3.9);
\draw (-1,-4.1) -- (-2,-4.9);
\draw (-1,-4.1) -- (0,-4.9);
\draw (1,-4.1) -- (0,-4.9);
\draw (-2,-5.1) -- (-1,-5.9);
\draw (0,-5.1) -- (-1,-5.9);
\draw (0,-5.1) -- (1,-5.9);
\draw (-1,-6.1) -- (-2,-6.9);
\draw (-1,-6.1) -- (0,-6.9);
\draw (1,-6.1) -- (0,-6.9);
\draw (-2,-7.1) -- (-1,-7.9);
\draw (0,-7.1) -- (-1,-7.9);
\draw (-1,-8.1) -- (-1,-9.4);
\draw (-1,-9.6) -- (-1,-10.9);

\draw (0,-1.4) -- (-0.1,-1.3); 
\draw (0,-1.4) -- (0.1,-1.3); 
\draw (0,-2.9) -- (-0.1,-2.8); 
\draw (0,-2.9) -- (0.1,-2.8); 
\draw (-1,-3.9) -- (-0.9,-3.7); 
\draw (-1,-3.9) -- (-0.8,-3.9); 
\draw (1,-3.9) -- (0.8,-3.9); 
\draw (1,-3.9) -- (1,-3.7); 
\draw (-2,-4.9) -- (-1.9,-4.7); 
\draw (-2,-4.9) -- (-1.8,-4.9); 
\draw (0,-4.9) -- (-0.2,-4.9); 
\draw (0,-4.9) -- (0,-4.7); 
\draw (0,-4.9) -- (0.1,-4.7); 
\draw (0,-4.9) -- (0.2,-4.9); 
\draw (-1,-5.9) -- (-1.2,-5.9); 
\draw (-1,-5.9) -- (-1,-5.7); 
\draw (-1,-5.9) -- (-0.9,-5.7); 
\draw (-1,-5.9) -- (-0.8,-5.9); 
\draw (1,-5.9) -- (0.8,-5.9); 
\draw (1,-5.9) -- (1,-5.7); 
\draw (-2,-6.9) -- (-1.9,-6.7); 
\draw (-2,-6.9) -- (-1.8,-6.9); 
\draw (0,-6.9) -- (-0.2,-6.9); 
\draw (0,-6.9) -- (0,-6.7); 
\draw (0,-6.9) -- (0.1,-6.7); 
\draw (0,-6.9) -- (0.2,-6.9); 
\draw (-1,-7.9) -- (-1.2,-7.9); 
\draw (-1,-7.9) -- (-1,-7.7); 
\draw (-1,-7.9) -- (-0.9,-7.7); 
\draw (-1,-7.9) -- (-0.8,-7.9); 
\draw (-1,-9.4) -- (-1.1,-9.3); 
\draw (-1,-9.4) -- (-0.9,-9.3); 
\draw (-1,-10.9) -- (-1.1,-10.8); 
\draw (-1,-10.9) -- (-0.9,-10.8); 

\end{tikzpicture} 
\caption{Diagram of $\Delta_n^+$}\label{diagram}
\end{center}
\end{figure}

 In the Figure \ref{diagram}, the vertices represent the roots in $\Delta_n^+$. Two roots $\beta, \gamma \in \Delta_n^+$ are joined by an edge with an arrow in the direction of 
$\gamma$ if $\gamma = \beta + \phi$ for some simple root $\phi \in \Delta_\frak{k}^+$. In this case, the simple root $\phi$ is given on one side of the edge. 

Let $\frak{b}_\frak{k}$ be the Borel subalgebra $\frak{h} \oplus \sum_{\alpha \in \Delta_\frak{k}^+} \frak{g}^\alpha$ of $\frak{k},$ $\mathcal{B}$ be the set of all Borel subalgebras 
of $\frak{g}$ which contain $\frak{b}_\frak{k},$ and $\mathcal{Q}$ be the set of all parabolic subalgebras of $\frak{g}$ containing a Borel subalgebra of $\mathcal{B}.$ 
Now a $\theta$-stable parabolic subalgebra of $\frak{g}_0$ is a parabolic subalgebra of $\frak{g},$ and a parabolic subalgebra of $\frak{g}$ is $\theta$-stable {\it if and only if} it contains a maximal 
abelian subspace of $\frak{k}_0.$ If $\frak{q}$ is a $\theta$-stable parabolic subalgebra of $\frak{g}_0,$ so is $Ad(k)(\frak{q}),$ and $A_\frak{q}$ is unitarily equivalent to $A_{Ad(k)(\frak{q})}$ for all $k \in K.$ 
Thus it is sufficient to consider parabolic subalgebras of $\frak{g}$ which are in $\mathcal{Q}.$  
Let $\frak{q}$ be a parabolic subalgebra of $\frak{g}.$ Then $\frak{q} \in \mathcal{Q}$ {\it if and only if} there exists a positive root system 
$\Delta_\frak{q}^+$ of $\Delta$ containing $\Delta_\frak{k}^+$ and a subset $\Gamma$ of $\Phi_\frak{q}$, the set of all simple roots in $\Delta_\frak{q}^+$, such that 
$\frak{q}= \frak{h} \oplus \sum_{\alpha \in \Delta, n_\phi(\alpha) \ge 0 \textrm{ for all }\alpha \in \Gamma} \frak{g}^\alpha,$ where $\alpha = \sum_{\phi \in \Phi_\frak{q}}n_\phi(\alpha) \phi 
\in \Delta.$ So we will determine all positive root systems of $\Delta$ containing $\Delta_\frak{k}^+.$ The number of such positive root systems is $|W_{\frak{f}_4}|/|W_{\frak{sp}(3)\oplus \frak{su}(2)}|
=1152/(2^3\times3! \times 2!)=12.$ 
%See \cite[Lemma 2.6]{paul2}. 

Let $\Phi$ be a basis of the root system $\Delta$ such that the corresponding positive root system $P(\Phi)$ contains $\Delta_\frak{k}^+.$ We will get all such $\Phi$ and their 
corresponding $P(\Phi)$ from the Figure \ref{diagram}, and these are listed in the Table \ref{+system}.

\newpage
\begin{table}[!h]
\caption{Positive root systems and corresponding simple systems :}\label{+system}
\begin{tabular}{||c|c|c||}
\hline
\  & $P(\Phi)$ & $\Phi$ \\
\hline
\hline
$1.$ & $\Delta_\frak{k}^+ \cup \{\beta \in \Delta_n^+ : \beta \ge \phi_1\}$ & $\{\phi_1, \phi_2, \phi_3, \phi_4 \}$ \\
\hline 
$2.$ & \begin{tabular}{c} $\Delta_\frak{k}^+ \cup \{\beta \in \Delta_n^+ : \beta \ge \phi_1+\phi_2\}$ \\ $\cup (-\{\beta \in \Delta_n^+ : \beta \le \phi_1\})$ \end{tabular} & 
$\{\phi_1+\phi_2,-\phi_1,\phi_3, \phi_4\}$ \\
\hline
$3.$ & \begin{tabular}{c} $\Delta_\frak{k}^+ \cup \{\beta \in \Delta_n^+ : \beta \ge \phi_1+\phi_2+\phi_3\}$ \\ $\cup (-\{\beta \in \Delta_n^+ : \beta \le \phi_1+\phi_2\})$ \end{tabular} & 
$\{\phi_1+\phi_2+\phi_3,-(\phi_1+\phi_2),\phi_2, \phi_4\}$ \\
\hline
$4.$ & \begin{tabular}{c} $\Delta_\frak{k}^+ \cup \{\beta \in \Delta_n^+ : \beta \ge \phi_1+\phi_2+2\phi_3 $\\$\textrm{ or } \phi_1+\phi_2+\phi_3+\phi_4\}$ \\ $\cup (-\{\beta \in \Delta_n^+ : \beta \le \phi_1+\phi_2+\phi_3\})$ \end{tabular} & 
\begin{tabular}{c}$\{\phi_1+\phi_2+2\phi_3,\phi_1+\phi_2+\phi_3+\phi_4,$\\$-(\phi_1+\phi_2+\phi_3),\phi_2\}$ \end{tabular} \\
\hline
$5.$ & \begin{tabular}{c} $\Delta_\frak{k}^+ \cup \{\beta \in \Delta_n^+ : \beta \ge \phi_1+\phi_2+2\phi_3\}$ \\ $\cup (-\{\beta \in \Delta_n^+ : \beta \le \phi_1+\phi_2+\phi_3+\phi_4\})$ \end{tabular} & 
$\{\phi_1+\phi_2+2\phi_3,-(\phi_1+\phi_2+\phi_3+\phi_4),\phi_2,\phi_4\}$ \\
\hline
$6.$ & \begin{tabular}{c} $\Delta_\frak{k}^+ \cup \{\beta \in \Delta_n^+ : \beta \ge \phi_1+2\phi_2+2\phi_3 $\\$\textrm{ or } \phi_1+\phi_2+\phi_3+\phi_4\}$ \\ $\cup (-\{\beta \in \Delta_n^+ : \beta \le \phi_1+\phi_2+2\phi_3\})$ \end{tabular} & 
\begin{tabular}{c}$\{\phi_1+2\phi_2+2\phi_3,\phi_1+\phi_2+\phi_3+\phi_4,$\\$-(\phi_1+\phi_2+2\phi_3),\phi_3\}$ \end{tabular} \\
\hline
$7.$ & \begin{tabular}{c} $\Delta_\frak{k}^+ \cup \{\beta \in \Delta_n^+ : \beta \ge \phi_1+\phi_2+\phi_3+\phi_4\}$ \\ $\cup (-\{\beta \in \Delta_n^+ : \beta \le \phi_1+2\phi_2+2\phi_3\})$ \end{tabular} & 
$\{\phi_1+\phi_2+\phi_3+\phi_4,-(\phi_1+2\phi_2+2\phi_3),\phi_2,\phi_3\}$ \\
\hline
$8.$ & \begin{tabular}{c} $\Delta_\frak{k}^+ \cup \{\beta \in \Delta_n^+ : \beta \ge \phi_1+2\phi_2+2\phi_3 $\\$\textrm{ or } \phi_1+\phi_2+2\phi_3+\phi_4\}\cup 
(-\{\beta \in \Delta_n^+ : $\\$\beta \le \phi_1+\phi_2+2\phi_3 \textrm{ or }\phi_1+\phi_2+\phi_3+\phi_4\})$ \end{tabular} & 
\begin{tabular}{c}$\{\phi_1+2\phi_2+2\phi_3,\phi_1+\phi_2+2\phi_3+\phi_4,$\\$-(\phi_1+\phi_2+2\phi_3),-(\phi_1+\phi_2+\phi_3+\phi_4)\}$ \end{tabular} \\
\hline
$9.$ & \begin{tabular}{c} $\Delta_\frak{k}^+ \cup \{\beta \in \Delta_n^+ : \beta \ge \phi_1+2\phi_2+2\phi_3 $\\$\textrm{ or } \phi_1+\phi_2+2\phi_3+2\phi_4\}$ \\ $\cup (-\{\beta \in \Delta_n^+ : \beta \le \phi_1+\phi_2+2\phi_3+\phi_4\})$ \end{tabular} & 
\begin{tabular}{c}$\{\phi_1+2\phi_2+2\phi_3,\phi_1+\phi_2+2\phi_3+2\phi_4,$\\$-(\phi_1+\phi_2+2\phi_3+\phi_4),\phi_3\}$ \end{tabular} \\
\hline
$10.$ & \begin{tabular}{c} $\Delta_\frak{k}^+ \cup \{\beta \in \Delta_n^+ : \beta \ge \phi_1+\phi_2+2\phi_3+\phi_4\}$ \\ $\cup (-\{\beta \in \Delta_n^+ : \beta \le \phi_1+2\phi_2+2\phi_3 $\\$\textrm{ or }\phi_1+\phi_2+\phi_3+\phi_4\})$ \end{tabular} & 
\begin{tabular}{c}$\{\phi_1+\phi_2+2\phi_3+\phi_4,-(\phi_1+2\phi_2+2\phi_3),$\\$-(\phi_1+\phi_2+\phi_3+\phi_4),\phi_2\}$ \end{tabular} \\
\hline
$11.$ & \begin{tabular}{c} $\Delta_\frak{k}^+ \cup \{\beta \in \Delta_n^+ : \beta \ge \phi_1+2\phi_2+2\phi_3\}$ \\ $\cup (-\{\beta \in \Delta_n^+ : \beta \le \phi_1+\phi_2+2\phi_3+2\phi_4\})$ \end{tabular} & 
$\{\delta,-(\phi_1+\phi_2+2\phi_3+2\phi_4),\phi_3,\phi_4\}$ \\
\hline
$12.$ & \begin{tabular}{c} $\Delta_\frak{k}^+ \cup \{\beta \in \Delta_n^+ : \beta \ge \phi_1+2\phi_2+2\phi_3+\phi_4 $\\$\textrm{ or } \phi_1+\phi_2+2\phi_3+2\phi_4\}\cup 
(-\{\beta \in \Delta_n^+ : $\\$\beta \le \phi_1+2\phi_2+2\phi_3 \textrm{ or }\phi_1+\phi_2+2\phi_3+\phi_4\})$ \end{tabular} & 
\begin{tabular}{c}$\{\delta,-(\phi_1+2\phi_2+2\phi_3),$\\$-(\phi_1+\phi_2+2\phi_3+\phi_4),\phi_3\}$ \end{tabular} \\
\hline
\end{tabular} 
\end{table} 

\begin{remark} 
Since $(\phi_1+\phi_2+2\phi_3+\phi_4)+(\phi_1+2\phi_2+2\phi_3+\phi_4)= \delta= (\phi_1+2\phi_2+2\phi_3)+(\phi_1+\phi_2+2\phi_3+2\phi_4),$ so if $-(\phi_1+\phi_2+2\phi_3+\phi_4)$ and $-(\phi_1+2\phi_2+2\phi_3+\phi_4),$ 
or $-(\phi_1+2\phi_2+2\phi_3)$ and $-(\phi_1+\phi_2+2\phi_3+2\phi_4)$ lie in some $P(\Phi),$ then $-\delta \in P(\Phi),$ and hence $P(\Phi)$ cannot contain $\Delta_\frak{k}^+.$ 
\end{remark}

Let $\Phi$ be a basis of the root system $\Delta$ such that the corresponding positive root system $P(\Phi)$ contains $\Delta_\frak{k}^+,$ and 
$\frak{b}_\Phi = \frak{h} \oplus \sum_{\alpha \in P(\Phi)} \frak{g}^\alpha$ be the associated Borel subalgebra of $\frak{g}.$ If $\frak{q} \in \mathcal{Q}$ 
is a parabolic subalgebra of $\frak{g}$ which contains the Borel subalgebra $\frak{b}_\Phi,$ then there exists a subset $\Gamma$ of $\Phi$ such that 
$\frak{q}= \frak{l} \oplus \frak{u},$
where $\frak{l} = \frak{h} \oplus \sum_{n_\phi (\alpha) =0 \textrm{ for all } \phi \in \Gamma} \frak{g}^\alpha$ is the Levi subalgebra of $\frak{q},$ and 
$\frak{u} = \sum_{n_\phi (\alpha) >0 \textrm{ for some } \phi \in \Gamma} \frak{g}^\alpha$ is the nilradical of $\frak{q}.$ The Levi subalgebra $\frak{l}$ 
is the direct sum of an $|\Gamma|$-dimensional centre and a semisimple Lie algebra whose Dynkin diagram is the subdiagram of the dynkin diagram 
of $\frak{g}$ consisting of the vertices $\Phi \setminus \Gamma.$ Write $\Delta(\frak{u} \cap \frak{p})=\{\beta \in \Delta_n: \frak{g}^\beta \subset \frak{u}\}.$ 
For $\frak{q}, \frak{q}' \in \mathcal{Q}, A_\frak{q}$ is unitarily equivalent to $A_{\frak{q}'}$ {\it if and only if} $\Delta(\frak{u} \cap \frak{p})=
\Delta(\frak{u'} \cap \frak{p})$ \cite{riba}. So we will define an equivalence relation on $\mathcal{Q}$ by $\frak{q} \sim \frak{q}'$ if 
$\Delta(\frak{u} \cap \frak{p})=\Delta(\frak{u}' \cap \frak{p}).$ For each basis $\Phi$ of the root system $\Delta$ such that the corresponding positive root system 
$P(\Phi)$ contains $\Delta_\frak{k}^+,$ we will determine the parabolic subalgebras $\frak{q} \in \mathcal{Q}$ which contain the Borel subalgebra $\frak{b}_\Phi,$ and their 
corresponding $\Delta(\frak{u} \cap \frak{p}).$ For an equivalence class of parabolic subalgebras of $\mathcal{Q},$ we will describe only one representative of the 
equivalence class. 

\noindent
{\bf $1. \Phi=\{\nu, \phi_2, \phi_3, \phi_4\}$, where $\nu = \phi_1$:} 
\begin{table}[!h]
 \\
\hline
\end{tabular} 
\end{table}

For $\frak{q} \in \mathcal{Q},$ let $R(\frak{q})= \textrm{dim}_\mathbb{C}(\frak{u} \cap \frak{p}),$ and $Y_\frak{q}$ be the compact Riemannian globally symmetric 
space of type $(\frak{l}_0^*, \frak{l}_0 \cap \frak{k}_0),$ where $\frak{u}$ is the nilradical of $\frak{q},$ $\frak{l}$ is the Levi subalgebra of $\frak{q}, \frak{l}_0 = \frak{l} \cap \frak{g}_0$  
is a real form of $\frak{l},$ and $\frak{l}_0^*$ is the compact real form of $\frak{l}$ dual to the pair $(\frak{l}_0, \frak{l}_0 \cap \frak{k}_0).$ Since $\delta \in \Delta_\frak{k}^+ \subset 
P(\Phi), \delta = \sum_{\phi \in \Phi}d_\phi \phi, d_\phi \in \mathbb{N} \cup \{0\},$ for all $\phi \in \Phi.$ So $\frak{g}^\delta \subset \frak{l}$ {\it iff} $\Gamma \subset 
\{\phi \in \Phi : d_\phi =0 \}.$ If  $\frak{g}^\delta \subset \frak{l},$ then $\delta \in \Delta([\frak{l},\frak{l}], [\frak{l},\frak{l}] \cap \frak{h})$ is a compact root, and $\delta$ is the  
highest root of a simple ideal of $[\frak{l}, \frak{l}]$ relative to the positive root system $P(\Phi) \cap \Delta([\frak{l},\frak{l}], [\frak{l},\frak{l}] \cap \frak{h}).$ See Table \ref{f-table}. 
If $P_\frak{q}(t)$ denotes the Poincar\'{e} polynomial of the cohomologies $ H^* (\frak{g}, K; A_{\frak{q}, K})$, then we have 
\[ P_\frak{q}(t) = t^{R(\frak{q})} P(Y_\frak{q} , t). \]
In the Table \ref{f-table}, we have listed $\Phi, \Gamma, \Delta(\frak{u} \cap \frak{p})$ for each $\frak{q} \in \mathcal{Q},$ and we have determined 
$Y_\frak{q}$ for each such $\frak{q}.$ We can see that $Y_\frak{q}$ is either singleton, or $\frac{SU(k)}{S(U(1)\times U(k-1))}(k=2,3),$ or $\frac{SO(2k+1)}{SO(2)\times SO(2k-1)}
(k=2,3),$ or $\frac{SO(7)}{SO(4)\times SO(3)},$ or $\frac{Sp(k)}{Sp(1)\times Sp(k-1)}(k =2,3),$ or $\frac{Sp(3)}{U(3)},$ or $\frac{F_4}{Sp(3) \times SU(2)}$($F_4$ is a 
compact connected Lie group with Lie algebra the compact real form of $\frak{f}_4$), or $\frac{SU(2)}{SO(2)} \times \frac{SU(3)}{S(U(1)\times U(2))},$ or 
$\frac{SU(2)}{SO(2)} \times \frac{SU(2)}{SO(2)}.$

We have $\frac{SU(2)}{SO(2)} \cong \frac{SU(2)}{S(U(1)\times U(1))},$ and 
$P(\frac{SU(k)}{S(U(1)\times U(k-1))},t)=1+t^2+t^4+\cdots+t^{2k-2} \textrm{ for }k \ge 2$ \cite[Table I, Ch. XI]{ghv}. So $P(\frac{SU(2)}{SO(2)} \times \frac{SU(3)}{S(U(1)\times U(2))} , t) = 
(1+t^2)(1+t^2+t^4)=1+2t^2+2t^4+t^6,$ and $P(\frac{SU(2)}{SO(2)} \times \frac{SU(2)}{SO(2)}, t) = (1+t^2)(1+t^2)=1+2t^2+t^4.$
Also $P(\frac{SO(2k+1)}{SO(2)\times SO(2k-1)},t)=1+t^2+t^4+\cdots+t^{4k-2} \textrm{ for all }k \ge 1,$ and $P(\frac{SO(7)}{SO(4)\times SO(3)} , t) = 1+2t^4+2t^8+t^{12}$ 
\cite[Table II, Ch. XI]{ghv}. 
Again $P(\frac{Sp(k)}{Sp(1)\times Sp(k-1)},t)=1+t^4+\cdots+t^{4k-4} \textrm{ for all }k \ge 2,$ and $P(\frac{Sp(3)}{U(3)}, t) = 1+t^2+t^4+2t^6+t^8+t^{10}+t^{12}$ \cite[Table IV, Ch. XI]{ghv}, 
$P(\frac{F_4}{Sp(3) \times SU(2)} , t) = 1+t^4+2t^8+2t^{12}+2t^{16}+2t^{20}+t^{24}+t^{28}$ \cite{pic}, and $P(\textrm{singleton},t)=1.$

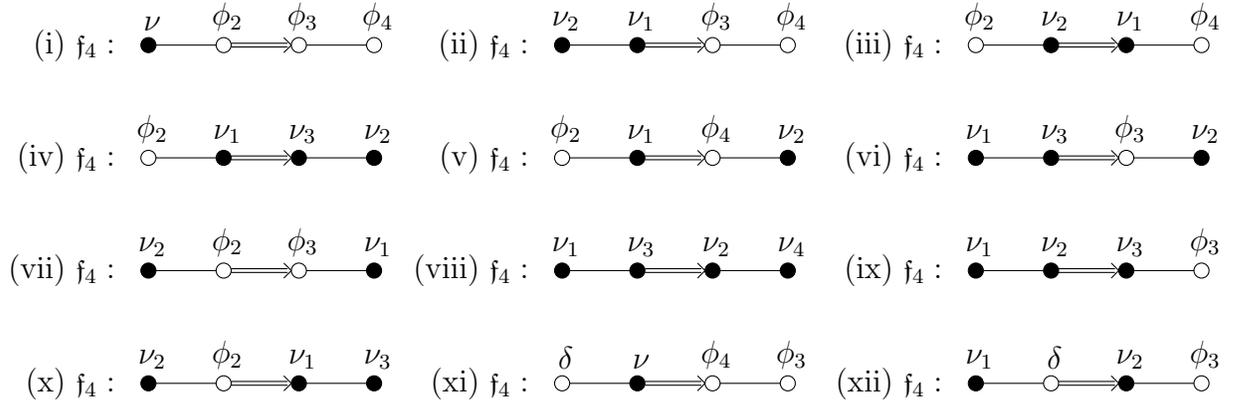
\begin{landscape} 
\begin{table} 
\caption{Poincar\'{e} Polynomial of The Irreducible Unitary Representations $A_\frak{q}$ of a connected Lie group of type $\frak{f}_{4(4)}:$}\label{f-table}

\end{center}

So the number of unitary equivalence classes of irreducible unitary representations with non-zero $(\frak{g},K)$-cohomology of a connected Lie group $G$ with Lie algebra 
$\frak{f}_{4(4)}$ is $46.$ From the Theorem \ref{th1}, the discrete series representations of $G$ with trivial infinitesimal character are 
given in the Table \ref{f-table} as 1.(b), 2.(b), 3.(b), 4.(e), 5.(b), 6.(e), 7.(b), 
8.(d), 9.(d), 10.(c), 11.(a), and 12.(a). Among these, the representation 1.(b) is the only Borel-de Siebenthal discrete series representation 
with trivial infinitesimal character, which is expected from \cite[Theorem 1.2]{paul2}.

\section*{acknowledgement}

The author acknowledges the financial support from the Department of Science and Technology (DST), Govt. of India under the Scheme 
``Fund for Improvement of S\&T Infrastructure (FIST)" [File No. SR/FST/MS-I/2019/41].

\end{document}